\documentclass[12pt]{article}
\usepackage{mathrsfs}
\usepackage{amsfonts}
\usepackage{amssymb,amsmath}
\usepackage{cite}
\usepackage{cases}
\usepackage[usenames]{color}

\def\UU{{\frak U}}

\def\Der{{\rm Der}}
\def\Inn{{\rm Inn}}
\def\Ker{{\rm Ker}}
\def\Im{{\rm Im}}

\def\cl{\centerline}
\def\rar{\rightarrow}

\def\vs{\vspace*}
\def\ni{\noindent}
\def\VV{\mathcal {V}}
\def\Z{\mathbb{Z}}
\def\C{\mathbb{C}}
\def\QED{\hfill$\Box$}

\def\D{\Delta}
\def\LL{{\cal L}}
\def\sl{\frak{sl}_2(C_q[x,y])}

\def\Cq{C_q[x,y]}
\def\sm{\setminus}
\def\d{\delta}
\def\H{\frak{H}}
\def\cq{{\cal C}_q}
\def\bz{{\bf 0}}
\def\bc{{\bf c}}
\def\bk{{\bf k}}

\def\bm{{\bf m}}
\def\DD{{\cal D}}
\def\dir{{\mathscr{D}}}
\def\bn{{\bf n}}
\def\EE{{\cal E}}
\def\FF{{\cal F}}
\def\GG{{\cal G}}
\def\HH{{\cal H}}
\def\WL{\widetilde{\cal{L}}}
\def\NN{{\cal N}}
\def\ek{e_{\bf k}}
\def\em{e_{\bf m}}

\def\fk{f_{\bf k}}
\def\fm{f_{\bf m}}

\def\ep{\epsilon}

\def\eki{\epsilon_i({\bf k})}

\def\eko{\epsilon_1({\bf k})}
\def\ekt{\epsilon_2({\bf k})}
\def\emi{\epsilon_i({\bf m})}
\def\emj{\epsilon_j({\bf m})}
\def\emo{\epsilon_1({\bf m})}
\def\emt{\epsilon_2({\bf m})}

\numberwithin{equation}{section}
\newtheorem{theo}{Theorem}[section]
\newtheorem{defi}[theo]{Definition}
\newtheorem{coro}[theo]{Corollary}
\newtheorem{lemm}[theo]{Lemma}
\newtheorem{prop}[theo]{Proposition}
\newtheorem{clai}{Claim}

\newtheorem{subclai}{Subclaim}

\begin{document}
\cl{{\large\bf Lie bialgebra structures on the Lie algebra $\widetilde{\frak{sl}_2(C_q[x,y])}$}
\footnote{\!\!\!\!Supported by NSF grants (No 10825101, 11101056) and the China Postdoctoral Science Foundation Grant (No 201003326)}}\vs{6pt}

\cl{Ying Xu$^{1)}$, Junbo Li$^{2)}$, Wei Wang$^{3)}$}\vs{6pt}

\cl{\small $^{1)}$Wu Wen-Tsun Key Laboratory of \vs{-2pt}Mathematics,}
\cl{\small  University of Science and
Technology of China, Hefei 230026, China}

\cl{\small $^{2)}$School\! of\! Mathematics\! and\! Statistics, Changshu Institute of Technology, Changshu 215500, China}

\cl{\small $^{3)}$School of Mathematics and Computer Science,
Ningxia University, Yinchuan 750021, China}\vs{6pt}

{\small
\parskip .005 truein
\baselineskip 3pt \lineskip 3pt

\noindent{{\bf Abstract.}\ \ In the present paper we shall
investigate the Lie bialgebra structures on the Lie algebra $\widetilde{\frak{sl}_2(C_q[x,y])}$, which are shown to be triangular coboundary.\vs{5pt}

}

}
\parskip .001 truein\baselineskip 6pt \lineskip 6pt

\vs{18pt}

\cl{\bf\S1. \
Introduction}\setcounter{section}{1}\setcounter{equation}{0}

\vs{8pt}

In order to find more about the solutions of the Yang-Baxter quantum equations, Drinfeld firstly introduced the notion of Lie bialgebras in 1983 (see \cite{D1}). Since then this issue has caused wide public concern (see \cite{D2},\cite{G},\cite{Mi}). Witt and Virasoro type Lie bialgebras were introduced in \cite{T}, which were further classified in \cite{NT}. The generalized case was considered in \cite{SS}. Lie bialgebra structures on generalized Virasoro-like and Block Lie algebras were investigated in \cite{WSS} and \cite{LSX} respectively. The same problem on the $q$-analog Virasoro-like algebra was settled in \cite{CS}. In the paper we shall concentrate on this problem on a more completed Lie algebra $\widetilde{\frak{sl}_2(C_q[x,y])}$, which is closely related to the $q$-analog Virasoro-like algebra.

Fix a nonzero complex number $q$ which is not a root of unity. Let $\cq$ be the $\C$-algebra $\Cq$ defined by generators $x$, $y$ and relations $yx=qxy$. Firstly, we introduce the Lie algebra $\LL=\sl$. Set $\bz=(0,0)$. For ${\bf k}=(k_1,k_2)\in\Z^2_{\geq 0}$, ${\bf m}=(m_1,m_2)\in\Z^2_{\geq 0}\!\sm\!\{\bz\}$ and $i=1,2$, define the following elements of $\LL$
\begin{eqnarray*}
\ek=E_{12}x^{k_1}y^{k_2},\ \ \ \fk=E_{21}x^{k_1}y^{k_2},\ \ \ h=E_{11}-E_{22},\ \ \ \emi=E_{ii}x^{m_1}y^{m_2},
\end{eqnarray*}
which form a basis of $\LL$ and satisfy the following brackets
\begin{eqnarray*}
&&[\eko,\fm]=-q^{k_1m_2}f_{\bk+\bm},\ \ \ [\ekt,\fm]=q^{k_2m_1}f_{\bk+\bm},\\
&&[\eko,\em]=q^{k_2m_1}e_{\bk+\bm},\ \ \ \ \ \,[\ekt,\em]=-q^{k_1m_2}e_{\bk+\bm},\\
&&[\ek,\fm]=q^{k_2m_1}\ep_1(\bk+\bm)-q^{m_2k_1}\ep_2(\bk+\bm),\\
&&[\eki,\emj]=\d_{ij}(q^{k_2m_1}-q^{m_2k_1})\ep_i(\bk+\bm),\\
&&[\ek,\em]=[\fk,\fm]=[h,\emi]=0,\ \ [h,\ek]=2\ek,\ \ [h,\fk]=-2\fk.
\end{eqnarray*}
This Lie algebra appeared during the study of some extended affine Lie algebras in \cite{B} and while taking $q\!=\!1$ limit of the quantum toroidal algebras in \cite{V}. Representation and central extensions of the loop algebra with the algebra of Laurent polynomials replaced by a quantum torus \cite{MY} attract a great many specialists and scholars (see\cite{J}, \cite{M}, \cite{BJ}, \cite{G1} and \cite{G2}). In this paper, we investigate the Lie bialgebra structures on $\widetilde{\sl}$ with replacing Laurent polynomials by polynomials.

For convenience we introduce the following notations
\begin{eqnarray*}
\EE_{k_1,k_2}=\ek,\ \ \ \FF_{k_1,k_2}=\fk,\ \ \ \GG_{m_1,m_2}=\emo,\ \ \ \HH_{m_1,m_2}=\emt,\ \ \ {\DD}=h.
\end{eqnarray*}
Then the relations can be rewritten as follows under the new basis notations
\begin{eqnarray*}
&&\!\!\!\!\!\!\!\![\EE_{k_1,k_2},\EE_{l_1,l_2}]=[\FF_{k_1,k_2},\FF_{l_1,l_2}]=0,\ \ \ [\GG_{m_1,m_2},\HH_{n_1,n_2}]=0,\\
&&\!\!\!\!\!\!\!\![\DD,\EE_{k_1,k_2}]=2\EE_{k_1,k_2},\ \ [\DD,\FF_{k_1,k_2}]=-2\FF_{k_1,k_2},\ \ [\DD,\GG_{m_1,m_2}]=[\DD,\HH_{m_1,m_2}]
=0,\\
&&\!\!\!\!\!\!\!\![\GG_{m_1,m_2},\FF_{k_1,k_2}]=-q^{m_1k_2}\FF_{m_1+k_1,m_2+k_2},\ \ \ [\HH_{m_1,m_2},\FF_{k_1,k_2}]=q^{m_2k_1}\FF_{m_1+k_1,m_2+k_2},\\
&&\!\!\!\!\!\!\!\![\GG_{m_1,m_2},\EE_{k_1,k_2}]=q^{m_2k_1}\EE_{m_1+k_1,m_2+k_2},\ \ \ \ \ \ \,[\HH_{m_1,m_2},\EE_{k_1,k_2}]=-q^{m_1k_2}\EE_{m_1+k_1,m_2+k_2},\\
&&\!\!\!\!\!\!\!\![\GG_{m_1,m_2},\GG_{n_1,n_2}]=(q^{m_2n_1}-q^{n_2m_1})\GG_{m_1+n_1,m_2+n_2},\\
&&\!\!\!\!\!\!\!\![\HH_{m_1,m_2},\HH_{n_1,n_2}]=(q^{m_2n_1}-q^{n_2m_1})\HH_{m_1+n_1,m_2+n_2},\\
&&\!\!\!\!\!\!\!\![\EE_{k_1,k_2},\FF_{l_1,l_2}]=
\left\{\begin{array}{cc}
q^{k_2{l_1}}\GG_{k_1+{l_1},k_2+{l_2}}-q^{{l_2k_1}}\HH_{k_1+{l_1},k_2+{l_2}}&{\rm if}\ \,(k_1+{l_1},k_2+{l_2})\neq(0,0),\\[6pt]
q^{k_2{l_1}}\DD&{\rm if}\ \,(k_1+{l_1},k_2+{l_2})=(0,0).
\end{array}\right.\\
\end{eqnarray*}
Introduce two degree derivations $\DD_1$ and $\DD_2$ on $\LL$, i.e.,
\begin{eqnarray*}
&&[\DD_1,\DD_2]=[\DD_i,\DD]=0,\\
&&[\DD_i,x]=k_i x \ \ {\rm for}\ \ x\in\LL_{k_1,k_2},\ \ i=1,2.
\end{eqnarray*}
Then we arrive at the Lie algebra $\WL=\LL\oplus\C\DD_1\oplus\C\DD_2$ which we shall consider in this paper. It is easy to see that $\WL$ possesses the following triangular decomposition: $\WL=\WL_-\oplus{\H}\oplus\WL_+$ with
\begin{eqnarray*}
&&\WL_-=\bigoplus\C\FF_{k_1,k_2},\ \ \WL_+=\bigoplus\C\EE_{k_1,k_2},\ \ \H={\frak{H}}_0\oplus\GG\oplus\HH,
\end{eqnarray*}
where ${\frak{H}}_0=\C\DD\oplus\C\DD_1\oplus\C\DD_2$, $\GG=\bigoplus\limits_{(m_1,m_2)\neq(0,0)}\C\GG_{m_1,m_2}$, $\HH=\bigoplus\limits_{(m_1,m_2)\neq(0,0)}\C\HH_{m_1,m_2}$. One can find the maximal commutative subalgebra of $\WL$, denoted $\NN=\C\DD\oplus\C\DD_1\oplus\bigoplus\limits_{m\neq0}\C\GG_{0,m}\oplus\bigoplus\limits_{m\neq0}\C\HH_{0,m}$. And $\WL$ is $\Z\times\Z$-graded: $\WL=\bigoplus\limits_{(k_1,k_2)\in\Z_{\geq 0}^2}\WL_{k_1,k_2}$ with $\WL_{k_1,k_2}=\C\EE_{k_1,k_2}\oplus\C\FF_{k_1,k_2}\oplus\C\GG_{k_1,k_2}\oplus\C\HH_{k_1,k_2}$ for $(k_1,k_2)\neq(0,0)$ and $\WL_{0,0}=\C\EE_{0,0}\oplus\C\FF_{0,0}\oplus\C\DD\oplus\C\DD_1\oplus\C\DD_2$.

Recall the relevant knowledge on Lie bialgebras, which could be found in \cite{LSX}  or \cite{WSS}.
For any $\C$-vector space ${\frak L}$, denote $\xi$ the {\it cyclic map} of ${\frak L}\otimes{\frak L}\otimes {\frak L}$ with $ \xi (x_{1} \otimes x_{2} \otimes x_{3}) =x_{2} \otimes x_{3} \otimes x_{1}$ and $\tau$ the {\it twist map} of
${\frak L}\otimes{\frak L}$ with $\tau(x_1\otimes x_2)=x_2\otimes x_1$ for any $x_1,x_2,x_3\in{\frak L}$. A {\it Lie algebra} is a pair $({\frak L},\delta)$ of a vector space
${\frak L}$ and a linear map $\delta :{\frak L}\otimes{\frak L}\rar{\frak L}$ (the {\it bracket}) admitting
\begin{eqnarray*}
&&\Ker(\textbf{1}-\tau)\subset\Ker\,\delta,\ \ \ \delta \cdot (\textbf{1} \otimes
\delta)\cdot(\textbf{1}+\xi+\xi^{2})=0.
\end{eqnarray*}
A {\it Lie coalgebra} is a pair $({\frak L},\D)$ of a vector space ${\frak L}$
and a linear map $\D:{\frak L}\to{\frak L}\otimes{\frak L}$ (the {\it cobracket}) admitting
\begin{eqnarray}\label{cLie-s-s}
&&\Im\,\D\subset\Im(\textbf{1}-\tau),\ \ \ (1+\xi+\xi^{2})\cdot
(\textbf{1}\otimes\D)\cdot\D=0.
\end{eqnarray}
For a Lie algebra ${\frak L}$, we always use the symbol ``$\cdot$'' to stand for the {\it
diagonal adjoint action}
\begin{eqnarray*}
x\cdot(\mbox{$\sum\limits_{i}$}{a_{i}\otimes b_{i}})=
\mbox{$\sum\limits_{i}$}({[x,a_{i}]\otimes
b_{i}+a_{i}\otimes[x,b_{i}]}).
\end{eqnarray*}
\begin{defi}\rm
A {\it Lie bialgebra} is a triple $({\frak L},\delta,\D)$ admitting the
following conditions
\begin{eqnarray}
\!\!\!\!\!\!\!\!\!\!\!\!&& \mbox{$({\frak L}, \delta)$ is a
Lie algebra},\ \
\mbox{$({\frak L},\D)$ is a Lie coalgebra},\nonumber\\
\!\!\!\!\!\!\!\!\!\!\!\!&& \D  \delta (x,y) = x \cdot
\D y - y \cdot \D x,\ \ \forall\,\,x,y\in{\frak L}\ \ \mbox{(the compatibility
condition).}\label{bLie-d}
\end{eqnarray}
\end{defi}

Denote by $\UU$ the universal enveloping algebra of ${\frak L}$ and by $\textbf{1}$
the identity element of $\UU$. For any $r =\sum_{i} {a_{i} \otimes
b_{i}}\in{\frak L}\otimes{\frak L}$, define $r^{ij},\,c(r),\,i,j=1,2,3$ to be
elements of $\UU \otimes \UU \otimes \UU$
\begin{eqnarray*}
r^{12}\!\!\!&=&\!\!\!\mbox{$\sum \limits_{i}$}{a_{i}\otimes b_{i}
\otimes \textbf{1}}, \ \ r^{13}= \mbox{$\sum \limits_{i}$}{a_{i}\otimes \textbf{1}\otimes b_{i}},\\
r^{23}\!\!\!&=&\!\!\!\mbox{$\sum\limits_{i}$}{\textbf{1}\otimes
a_{i}\otimes b_{i}},\ \ \ \mbox{$c(r)=[r^{12},r^{13}]+[r^{12},r^{23}]
+[r^{13},r^{23}]$}.
\end{eqnarray*}
\begin{defi}
\label{def2} \rm (1) A {\it coboundary Lie bialgebra} is a $4$-tuple
$({\frak L}, \delta, \D,r),$ where $({\frak L},\delta,\D)$ is a Lie bialgebra
and $r \in \Im(1-\tau) \subset{\frak L}\otimes{\frak L}$ such that $\D=\D_r$ is
a {\it coboundary of $r$}, where $\D_r$ is defined by
\begin{eqnarray}
\D_r(x)=x\cdot r\mbox{\ \ for\ \ }x\in{\frak L}.\label{D-r}
\end{eqnarray}

 (2) A coboundary Lie bialgebra $({\frak L},\delta,\D,r)$
is called {\it triangular} if it satisfies the following {\it
classical Yang-Baxter Equation}
\begin{eqnarray}
c(r)=0.\label{CYBE}
\end{eqnarray}
\end{defi}

The main result of this paper can be formulated as follows.
\begin{theo}
Every Lie bialgebra on $\WL$ is triangular
coboundary.\label{main}
\end{theo}

\vs{16pt}

\cl{\bf\S2. \ Proof of the main result}\setcounter{section}{2}
\setcounter{theo}{0}\setcounter{equation}{0}

\vs{10pt}

The following lemma can be found in \cite{D1} or \cite{NT}.\vs{-6pt}
\begin{lemm}
\label{some} Let ${\frak L}$ be a Lie algebra and $r\in\Im(1-\tau)\subset
{\frak L}\otimes{\frak L}$.

{\rm(1)} The tripple $({\frak L},[\cdot,\cdot], \D_r)$ is
a Lie bialgebra if and only if $r$ satisfies $(\ref{CYBE})$.

{\rm(2)} For any $x\in{\frak L}$,
\begin{eqnarray}
(1+\xi+\xi^{2})\cdot(1\otimes\D)\cdot\D(x)=x\cdot
c(r).\label{add-c}
\end{eqnarray}
\end{lemm}

We also obtain the following lemma.\vs{-6pt}
\begin{lemm}\label{Legr}
Regard $\WL^{\otimes n}$ (the $n$ copies tensor product of $\WL$) as an
$\WL$-module under the adjoint diagonal action of $\WL$. If $r\in\WL^{\otimes n}$ satisfying $\WL_{\bz}\cdot r=0$, $\EE_{1,0}\cdot r=0$ and $\FF_{1,0}\cdot r=0$, then $r=0$. In particular, if $x\cdot r=0$ for all $x\in\WL$ and some $r\in\WL^{\otimes n}$, then $r=0$.
\end{lemm}
\noindent{\it Proof}.\ \ We can write $\WL^{\otimes n}$ as $\sum_{\bm}\WL^{\otimes n}_\bm$ with
\begin{eqnarray*}
\WL^{\otimes n}_{\bm}=\mbox{$\sum\limits_{\bm_1+\bm_2+\cdots+\bm_n=\bm}$}\WL_{\bm_1}\otimes\WL_{\bm_2}\otimes\cdots\otimes\WL_{\bm_n}.
\end{eqnarray*}
For any $r=\mbox{$\sum\limits_{r_{\bm}}$}\WL^{\otimes n}_{r_{\bm}}\in\WL^{\otimes n}$,
since $\DD_1\cdot r=\DD_2\cdot r=0$, we have $r_\bm=r_\bz$, i.e,
\begin{eqnarray}\label{eq-r}
r=\mbox{$\sum\limits_{\bm_1+\bm_2+\cdots+\bm_n=0}$}
r_{\bm_1,\bm_2,\cdots,\bm_n}L_{\bm_1}\otimes L_{\bm_2}\otimes\cdots\otimes L_{\bm_n},
\end{eqnarray}
where $r_{\bm_1,\bm_2,\cdots,\bm_n}\in\C$, $L_{\bm_j}\in\WL_{\bm_j}$ for $1\leq j\leq n$. Define a total order on $\Z^n$ by
\begin{eqnarray}\label{eq-order}
i<j\Longleftrightarrow\!\!\!\!\!\!&&|i|<|j|\ \ \ \ or\nonumber\\
\!\!\!\!\!\!&&|i|=|j|\ \ \mbox{\rm but there exists a $q$ such that }\ i_q <j_q\ \mbox{\rm and}\ i_p =j_p\ \mbox{\rm for}\ p<q,
\end{eqnarray}
where $i=(i_1,i_2,\cdots i_n)\in\Z^n$, $j=(i_1,j_2,\cdots j_n)\in\Z^n$, and $|i|=\sum_{p=1}^ni_p$, $|j|=\sum_{p=1}^ni_p$.

Choose the maximal summand appearing in \eqref{eq-r}, denoted by $(\bn_1,\bn_2,\cdots,\bn_n)$, under the convention of \eqref{eq-order}. Since $\EE_{1,0}\cdot r=0$, there is $L_{\bn_j}=\EE_{\bn_j}$ or $L_{\bn_j}=\GG_{k,0}+\HH_{k,0}$ for $1\leq j\leq n$ and some $k\in\Z$, otherwise we would obtain a higher summand.  Then, by $\FF_{1,0}\cdot r=0$, we obtain $L_{\bn_j}=0$, which gives $r=0$.\QED\vskip5pt

An element $r\in\Im(1-\tau)\subset\WL\otimes\WL$ is said to satisfy
the \textit{modified Yang-Baxter equation} if
\begin{eqnarray}
x\cdot c(r)=0,\ \,\forall\,\,x\in\WL.\label{MYBE}
\end{eqnarray}
According to Lemma \ref{Legr}, we immediately obtain
\begin{coro}\label{coro1}
Some $r\in\Im(1-\tau)\subset\WL\otimes\WL$ satisfies $(\ref{CYBE})$ if and only if it satisfies $(\ref{MYBE})$.
\end{coro}\vs{-6pt}

Regard $\VV=\WL\otimes\WL$ as a $\WL$-module under the adjoint
diagonal action. Denote by $\Der(\WL,\VV)$ the set of
\textit{derivations} $\dir:\WL\to\VV$, namely, $\dir$ is a linear map
satisfying
\begin{eqnarray}
\dir([x,y])=x\cdot \dir(y)-y\cdot \dir(x),\label{deriv}
\end{eqnarray}
and $\Inn(\WL,\VV)$ the set consisting of the derivations $v_{\rm
inn},\, v\in\VV$, where $v_{\rm inn}$ is the \textit{inner
derivation} defined by
\begin{eqnarray}\label{inn}
v_{\rm inn}:x\mapsto x\cdot v.
\end{eqnarray}
Then
\begin{eqnarray*}
H^1(\WL,\VV)\cong\Der(\WL,\VV)/\Inn(\WL,\VV),
\end{eqnarray*}
where $H^1(\WL,\VV)$ is the {\it first cohomology group} of the Lie
algebra $\WL$ with coefficients in the $\WL$-module $\VV$.

\begin{prop}
\label{proposition} $\Der(\WL,\VV)=\Inn(\WL,\VV)$, or equivalently,
$H^1(\WL,\VV)=0$.
\end{prop}
\noindent{\it Proof}.\ \ Note that $\VV=\WL\otimes\WL=\bigoplus\limits_{\bk\in\Z_{\geq 0}^2}\VV_\bk$ are $\Z^2$-graded with
\begin{eqnarray*}
\VV_\bk=\mbox{$\sum\limits_{\bm+\bn=\bk}$}\WL_\bm\otimes\WL_\bn\ \ \,\,\mbox{for}\ \,\,\bm,\bn\in\Z_{\geq 0}^2.
\end{eqnarray*}
A derivation $\dir\in\Der(\WL,\VV)$ is {\it homogeneous of degree
$\bk\in\Z^2$} if $\dir(\WL_\bn)\subset \VV_{\bk+\bn}$ for all $\bn\in\Z^2$.
Denote
\begin{eqnarray*}
\Der(\WL,\VV)_\bk=\{\dir\in\Der(\WL,\VV)\,|\,{\rm deg\,}\dir=\bk\}\mbox{ \ \ for \ } \bk\in\Z^2.
\end{eqnarray*}
Let $\dir\in\Der(\WL,\VV)$.
For any $\bk\in\Z^2$, we define the linear map
$\dir_\bk:\WL\rightarrow\VV$ as follows: For any $\mu\in\WL_\bn$ with
$\bn\in\Z^2$, write $\dir(\mu)=\sum_{\bm\in\Z^2}\mu_\bm$ with $\mu_\bm\in\VV_\bm$,
then we set $\dir_\bk(\mu)=\mu_{\bn+\bk}$. Obviously, $\dir_\bk\in
\Der(\WL,\VV)_\bk$ and we have
\begin{eqnarray}\label{summable}
\dir=\mbox{$\sum\limits_{\bk\in\Z^2}\dir_\bk$},
\end{eqnarray}
which holds in the sense that for every $\mu\in\WL$, only finitely
many $\dir_\bk(\mu)\neq 0,$ and $\dir(\mu)=\sum_{\bk\in\Z^2}\dir_\bk(\mu)$
(we call such a
sum in (\ref{summable}) {\it summable}).\\[8pt]

We shall prove this proposition by several claims.
\begin{clai}\label{clai1}
If $\bn\in\Z\times\Z\sm\{(0,0)\}$, then
$\dir_\bn\in\Inn(\WL,\VV)$.
\end{clai}
\noindent{\it Proof}.\ \ Denote $\Gamma=\{(k_{1},k_{2})\,|\,k_{i}\in\Z,\,i=1,2\}$ and $\texttt{T}={\rm Span}_\C\{\DD_1,\DD_2\}.$ Define the nondegenerate bilinear map form $\C^2\times\texttt{T} \longrightarrow \C,$ $\rho(\bc)=(\bc,\rho)=c_{1}\rho_{1}+c_{2}\rho_{2},$ for $\bc=(c_{1},c_{2})\in\C^2,$ $\rho=\rho_{1}\DD_1+\rho_{2}\DD_2\in\texttt{T}.$ By linear algebra, one can choose $\rho\in\texttt{T}$ with $\rho(\bn)\neq0$ for $\bn\in\Z\times\Z\sm\{(0,0)\}$. Denote $v=(\rho(\bn))^{-1}\dir_\bn(\rho)\in\WL_{\bn}.$ Then for any $x\in\WL_{\bk},\bk\in\Gamma,$ applying $\dir_\bn$ to $[\rho,x]=\rho(\bk)x,$ using $\dir_\bn(x)\in\VV_{\bn+\bk},$ we have $$\rho(\bn+\bk)\dir_\bn(x)-x\cdot\dir_\bn(\rho)=\rho\cdot\dir_\bn(x)-x\cdot\dir_\bn(\rho)=\rho(\bk)\dir_\bn(x),$$ i.e., $\dir_\bn(x)=v_{\rm inn}(x).$ Thus $\dir_\bn\in\Inn(\WL,\VV).$ \QED\vskip5pt
\begin{clai}
\label{clai2} $\dir_\bz(\DD_1)=\dir_\bz(\DD_2)=0.$
\end{clai}
\noindent{\it Proof}.\ \ Applying $\dir_\bz$ to $[\DD_1,x]=k_{1}x$ and $[\DD_2,x]=k_{2}x$ for $x\in\WL_{k_{1},k_{2}},$ we obtain that $x\cdot\dir_\bz(\DD_{1})=x\cdot\dir_\bz(\DD_{2})=0,$ Thus by Lemma \ref{Legr}, $\dir_\bz(\DD_1)=\dir_\bz(\DD_2)=0.$ \QED\vskip5pt
\begin{clai}
\label{clai3}  Replacing $\dir_\bz$ by $\dir_\bz-u_{\rm inn}$ for some $u\in
\VV_\bz$, one can suppose $\dir_\bz(\WL)=0,$ i.e., $\dir_\bz\in\Inn(\WL,\VV).$
\end{clai}
\vskip4pt
\par
The proof of this claim will be done by several Subclaims.
\begin{subclai}
\label{subclaim1}  By replacing $\dir_\bz$ by $\dir_\bz-u_{\rm inn}$ for some $u\in
\VV_\bz,$ one can suppose $\dir_\bz(\WL_{\bz})=0.$
\end{subclai}
\noindent{\it Proof}.\ \ It is well known that the first cohomology group of $sl_2(\C)$ vanishes on finite dimensional modules. And the subalgebra $\{\EE_{\bz},\FF_{\bz},\DD\}$ can be regarded as the Lie algebra $sl_2(\C)$. Then the subclaim follows from a very simple computation which shall be omitted here.\QED\vskip5pt

\begin{subclai}
\label{subclaim2} $\dir_{\bz}(u)=c\ u\cdot(\EE_{\bz}\otimes\FF_{\bz}+\FF_{\bz}\otimes\EE_{\bz}+1/2\DD\otimes\DD)$ for some $c\in\C,$ where $u\in\WL_{0,1}$ or $u\in\WL_{1,0}.$
\end{subclai}
\noindent{\it Proof}.\ \ For any $Y_{0,1}\in\WL_{0,1}$, denote $\dir_{\bz}(Y_{0,1})$ as
\begin{eqnarray}\begin{array}{llllllll}\label{E01}
&&y^{e00}_{e01}\EE_{0,0}\!\otimes\!\EE_{0,1}+y^{e01}_{e00}\EE_{0,1}\!\otimes\!\EE_{0,0}+y^{e00}_{f01}\EE_{0,0}\!\otimes\!\FF_{0,1}\!+\!y^{e01}_{f00}\EE_{0,1}\otimes\FF_{0,0} \\[6pt]
&&+y^{e00}_{g01}\EE_{0,0}\otimes\GG_{0,1}\!+\!y^{e00}_{h01}\EE_{0,0}\otimes\HH_{0,1}\!+\!y^{e01}_{d}\EE_{0,1}\otimes\DD\!+\!y^{e01}_{d_1}\EE_{0,1}\otimes\DD_{1} \\[6pt]
&&+y^{e01}_{d_2}\EE_{0,1}\!\otimes\!\DD_{2}\!+\!y^{f00}_{e01}\FF_{0,0}\!\otimes\!\EE_{0,1}\!\!+\!\!y^{f01}_{e00}\FF_{0,1}\!\otimes\!\EE_{0,0}\!+\!y^{f00}_{f01}\FF_{0,0}\otimes\FF_{0,1} \\[6pt]
&&+y^{f01}_{f00}\FF_{0,1}\!\otimes\!\FF_{0,0}\!+\!y^{f00}_{g01}\FF_{0,0}\!\otimes\!\GG_{0,1}\!+\!y^{f00}_{h01}\FF_{0,0}\!\otimes\!\!\HH_{0,1}\!+\!y^{f01}_d\FF_{0,1}\!\!\otimes\!\!\DD \\[6pt]
&&+y^{f01}_{d_1}\FF_{0,1}\!\otimes\!\DD_{1}+\!y^{f01}_{d_2}\FF_{0,1}\!\otimes\!\DD_{2}+y^{g01}_{e00}\GG_{0,1}\!\otimes\!\EE_{0,0}+\!y^{g01}_{f00}\GG_{0,1}\!\otimes\!\FF_{0,0} \\[6pt]
&&+y^{g01}_d\GG_{0,1}\otimes\DD\!+\!y^{g01}_{d_1}\GG_{0,1}\otimes\DD_{1}\!+y^{g01}_{d_2}\GG_{0,1}\!\otimes\DD_{2}\!+y^{h01}_{e00}\HH_{0,1}\otimes\EE_{0,0}\\[6pt]
&&+y^{h01}_{f00}\HH_{0,1}\!\otimes\!\FF_{0,0}\!+\!y^{h01}_d\HH_{0,1}\otimes\DD \!+\!y^{h01}_{d_1}\HH_{0,1}\!\otimes\!\DD_{1}\!+\!y^{h01}_{d_2}\HH_{0,1}\!\otimes\DD_{2}\\[6pt]
&&+y^d_{e01}\,\DD\otimes\EE_{0,1} +y^d_{f01}\,\DD\otimes\FF_{0,1}\,+y^{d}_{g01}\,\DD\otimes\GG_{0,1}+y^d_{h01}\,\DD\otimes\HH_{0,1} \\[6pt]
&&+y^{d_1}_{e01}\DD_{1}\otimes\EE_{0,1}+y^{d_1}_{f01}\DD_{1}\otimes\FF_{0,1}+y^{d_1}_{g01}\DD_{1}\otimes\!\GG_{0,1}\!+\!y^{d_1}_{h01}\DD_{1}\otimes\HH_{0,1} \\[6pt]
&&+y^{d_2}_{e01}\DD_{2}\otimes\EE_{0,1}+y^{d_2}_{f01}
\DD_{2}\otimes\FF_{0,1}+\!y^{d_2}_{g01}\DD_{2}\otimes\!\GG_{0,1}
\!+\!y^{d_2}_{h01}\DD_{2}\otimes\HH_{0,1}.
\end{array}\end{eqnarray}

For any $Y_{1,0}\in\WL_{1,0}$, denote $\dir_{\bz}(Y_{1,0})$ as
\begin{eqnarray}\begin{array}{llllllll}\label{E10}
&&y^{e00}_{e10}\EE_{0,0}\otimes\EE_{1,0}\!+\!y^{e10}_{e00}\EE_{1,0}\otimes\EE_{0,0}+y^{e00}_{f10}\EE_{0,0}\otimes\FF_{1,0}\!+\!y^{e10}_{f00}\EE_{1,0}\otimes\FF_{0,0} \\[6pt]
&&+y^{e00}_{g10}\EE_{0,0}\otimes\GG_{1,0}\!+\!y^{e00}_{h10}\EE_{0,0}\otimes\HH_{1,0}+y^{e10}_{d}\EE_{1,0}\otimes\DD+y^{e10}_{d_1}\EE_{1,0}\otimes\DD_{1} \\[6pt]
&&+y^{e10}_{d_2}\EE_{1,0}\!\otimes\DD_{2}\!+\!y^{f00}_{e10}\FF_{0,0}\otimes\EE_{1,0}\!+\!y^{f10}_{e00}\FF_{1,0}\otimes\EE_{0,0}\!+\!y^{f00}_{f10}\FF_{0,0}\!\otimes\FF_{1,0} \\[6pt]
&&+y^{f10}_{f00}\FF_{1,0}\!\otimes\!\FF_{0,0}\!+\!y^{f00}_{g10}\FF_{0,0}\!\otimes\GG_{1,0}\!+\!y^{f00}_{h10}\FF_{0,0}\otimes\HH_{1,0}\!+\!y^{f10}_{d}\FF_{1,0}\!\otimes\DD \\[6pt]
&&+y^{f10}_{d_1}\FF_{1,0}\otimes\DD_{1}\!+\!y^{f10}_{d_2}\FF_{1,0}\!\otimes\DD_{2}\!+\!y^{g10}_{e00}\GG_{1,0}\otimes\EE_{0,0}\!+\!y^{g10}_{f00}\GG_{1,0}\otimes\FF_{0,0} \\[6pt]
&&+y^{g10}_{d}\GG_{1,0}\otimes\DD+y^{g10}_{d_1}\GG_{1,0}\!\otimes\DD_{1}+y^{g10}_{d_2}\GG_{1,0}\otimes\DD_{2}+y^{h10}_{e00}\HH_{1,0}\otimes\EE_{0,0}\\[6pt]
&&+y^{h10}_{f00}\HH_{1,0}\otimes\FF_{0,0}\!+\!y^{h10}_{d}\HH_{1,0}\otimes\DD \!+\!y^{h10}_{d_1}\HH_{1,0}\otimes\DD_{1}\!+\!y^{h10}_{d_2}\HH_{1,0}\otimes\DD_{2}\\[6pt]
&&+y^{d}_{e10}\DD\,\otimes\,\EE_{1,0}\,+\,y^{d}_{f10}\DD\otimes\FF_{1,0}\,+\,y^{d}_{g10}\DD\otimes\GG_{1,0}\,+\,y^d_{h10}\DD\otimes\HH_{1,0} \\[6pt]
&&+y^{d_1}_{e10}\DD_{1}\otimes\EE_{1,0}+y^{d_1}_{f10}\DD_{1}\otimes\FF_{1,0}+y^{d_1}_{g10}\DD_{1}\otimes\GG_{1,0}+y^{d_1}_{h10}\DD_{1}\otimes\HH_{1,0} \\[6pt]
&&+y^{d_1}_{e10}\DD_{1}\otimes\EE_{1,0}+y^{d_2}_{f10}\DD_{2}\otimes\FF_{1,0}
+y^{d_2}_{g10}\DD_{2}\otimes\GG_{1,0}+y^{d_2}_{h10}\DD_{2}\otimes\HH_{1,0}. \end{array}\end{eqnarray}

Applying $\dir_{\bz}$ to $[\DD,\EE_{0,1}]\!=\!2\EE_{0,1},$ using Subclaim \ref{subclaim1} and expression \eqref{E01}, we can simplify $\dir_{\bz}(\EE_{0,1})$ as
\begin{eqnarray}\begin{array}{ccc}\label{E012}
&&e^{e00}_{g01}\EE_{0,0}\!\otimes\!\GG_{0,1}\!\!+\!\!e^{e00}_{h01}\EE_{0,0}\!\otimes\!\HH_{0,1}\!\!+\!\!e^{e01}_{d}\EE_{0,1}\!\otimes\!\DD\!\!+\!\!e^{e01}_{d_1}\EE_{0,1}\!\otimes\!\DD_{1}\!\!+\!\!e^{e01}_{d_2}\ \EE_{0,1}\otimes\DD_{2}\\[6pt]
&&+e^{g01}_{e00}\GG_{0,1}\!\otimes\!\EE_{0,0}\!\!+\!\!e^{h01}_{e00}\HH_{0,1}\!\otimes\!\EE_{0,0}\!\!+\!\!e^{d}_{e01}\DD\!\otimes\!\EE_{0,1}\!\!+\!\!e^{d_1}_{e01}\DD_{1}\!\otimes\!\EE_{0,1}\!\!+\!\!e^{d_2}_{e01}\DD_{2}\!\otimes\!\EE_{0,1}.
\end{array}\end{eqnarray}
Applying $\dir_{\bz}$ to $[\DD,\FF_{0,1}]\!=\!-2\FF_{0,1},$ using Subclaim \ref{subclaim1} and \eqref{E01}, we can simplify $\dir_{\bz}(\FF_{0,1})$ as
\begin{eqnarray}\begin{array}{ccc}\label{F012}
&&f^{f00}_{g01}\FF_{0,0}\!\otimes\!\GG_{0,1}\!\!+\!\!f^{f00}_{h01}\FF_{0,0}\!\otimes\!\HH_{0,1}\!\!+\!\!f^{f01}_{d}\FF_{0,1}\!\otimes\!\DD\!\!+\!\!f^{f01}_{d_1}\FF_{0,1}\!\otimes\!\DD_{1}\!\!+\!\!f^{f01}_{d_2}\FF_{0,1}\!\otimes\!\DD_{2} \\[6pt]
&&\!+\!f^{g01}_{f00}\GG_{0,1}\!\otimes\!\FF_{0,0}\!\!+\!\!f^{h01}_{f00}\HH_{0,1}\!\!\otimes\!\FF_{0,0}\!\!+\!\!f^{d}_{f01}\DD\!\!\otimes\!\FF_{0,1}\!\!+\!\!f^{d_1}_{f01}\DD_{1}\!\otimes\!\FF_{0,1}\!\!+\!\!f^{d_2}_{f01}\DD_{2}\!\otimes\!\FF_{0,1}.
\end{array}\end{eqnarray}
Applying $\dir_{\bz}$ to $[\DD,\GG_{0,1}]=0,$ using Subclaim \ref{subclaim1} and expression \eqref{E01}, we can simplify $\dir_{\bz}(\GG_{0,1})$ as
\begin{eqnarray}\begin{array}{cccc}\label{G012}
&&g^{e00}_{f01}\EE_{0,0}\!\otimes\!\FF_{0,1}\!+\!g^{e01}_{f00}\EE_{0,1}\!\otimes\!\FF_{0,0}\!+\!g^{f00}_{e01}\FF_{0,0}\otimes\EE_{0,1}\!+\!g^{f01}_{e00}\FF_{0,1}\otimes\EE_{0,0}\\[6pt]
&&+g^{g01}_{d}\GG_{0,1}\otimes\DD\!+\!g^{g01}_{d_1}\GG_{0,1}\otimes\DD_{1}+g^{g01}_{d_2}\GG_{0,1}\otimes\DD_{2}+g^{h01}_{d}\HH_{0,1}\otimes\DD \\[6pt]
&&+g^{h01}_{d_1}\HH_{0,1}\otimes\DD_{1}\!+\!g^{h01}_{d_2}\HH_{0,1}\otimes\DD_{2}\!+\!g^{d}_{g01}\DD\otimes\GG_{0,1}+g^{d}_{h01}\DD\otimes\HH_{0,1}\\[6pt]
&&+g^{d_1}_{g01}\DD_{1}\otimes\GG_{0,1}\!+\!g^{d_1}_{h01}\DD_{1}\otimes\HH_{0,1}\!+\!g^{d_2}_{g01}\DD_{2}\otimes\GG_{0,1}\!+\!g^{d_2}_{h01}\DD_{2}\!\otimes\!\HH_{0,1}. \end{array}\end{eqnarray}
Applying $\dir_{\bz}$ to $[\DD,\HH_{0,1}]=0,$ using subclaim \ref{subclaim1} and expression \eqref{E01}, we can simplify $\dir_{\bz}(\HH_{0,1})$ as
\begin{eqnarray}\begin{array}{cccc}\label{H012}
&&h^{e00}_{f01}\EE_{0,0}\otimes\FF_{0,1}\!+\!h^{e01}_{f00}\EE_{0,1}\otimes\FF_{0,0}\!+\!h^{f00}_{e01}\FF_{0,0}\otimes\EE_{0,1}\!+\!h^{f01}_{e00}\FF_{0,1}\otimes\EE_{0,0} \\[6pt]
&&+h^{g01}_d\GG_{0,1}\otimes\DD\,+\,h^{g01}_{d_1}\GG_{0,1}\otimes\DD_{1}+h^{g01}_{d_2}\GG_{0,1}\otimes\DD_{2}\,+h^{h01}_d\HH_{0,1}\otimes\DD \\[6pt]
&&+h^{h01}_{d_1}\HH_{0,1}\otimes\DD_{1}+h^{h01}_{d_2}\HH_{0,1}\otimes\DD_{2}+h^{d}_{g01}\DD\otimes\GG_{0,1}+h^{d}_{h01}\DD\otimes\HH_{0,1}\\[6pt]
&&+h^{d_1}_{g01}\DD_{1}\otimes\GG_{0,1}+h^{d_1}_{h01}\DD_{1}\otimes\HH_{0,1}\!+\!h^{d_2}_{g01}\DD_{2}\otimes\GG_{0,1}\!+\!h^{d_2}_{h01}\DD_{2}\otimes\HH_{0,1}.
\end{array}\end{eqnarray}
Applying $\dir_{\bz}$ to $[\DD,\EE_{1,0}]\!=\!2\EE_{1,0},$ using Subclaim \ref{subclaim1} and expression \eqref{E10}, we can simplify $\dir_{\bz}(\EE_{1,0})$ as
\begin{eqnarray}\begin{array}{cccc}\label{E102}
&&e^{e00}_{g10}\EE_{0,0}\!\otimes\!\GG_{1,0}\!+\!e^{e00}_{h10}\EE_{0,0}\!\otimes\!\HH_{1,0}\!+\!e^{e10}_d\EE_{1,0}\!\otimes\!\DD\!+\!e^{e10}_{d_1}\EE_{1,0}\!\otimes\!\DD_{1}\!+\!e^{e10}_{d_2}\EE_{1,0}\!\otimes\!\DD_{2} \\[6pt]
&&\!+\!e^{g10}_{e00}\GG_{1,0}\!\otimes\!\EE_{0,0}\!+\!e^{h10}_{e00}\HH_{1,0}\!\!\otimes\!\EE_{0,0}\!+\!e^{d}_{e10}\DD\!\!\otimes\!\EE_{1,0}\!+\!e^{d_1}_{e10}\DD_{1}\!\otimes\!\EE_{1,0}\!+\!e^{d_2}_{e10}\DD_{2}\!\otimes\!\EE_{1,0}.
\end{array}\end{eqnarray}
Applying $\dir_{\bz}$ to $[\DD,\FF_{1,0}]\!=\!-2\FF_{1,0},$ using Subclaim \ref{subclaim1} and expression \eqref{E10}, we can simplify $\dir_{\bz}(\FF_{1,0})$ as
\begin{eqnarray}\begin{array}{cccc}\label{F102}
&&f^{f00}_{g10}\FF_{0,0}\!\otimes\!\GG_{1,0}\!\!+\!\!f^{f00}_{h10}\FF_{0,0}\!\otimes\!\HH_{1,0}\!\!+\!\!f^{f10}_{d}\FF_{1,0}\!\otimes\!\DD\!\!+\!\!f^{f10}_{d_1}\FF_{1,0}\!\otimes\!\DD_{1}\!\!+\!\!f^{f10}_{d_2}\FF_{1,0}\!\otimes\!\DD_{2} \\[6pt]
&&\!+\!f^{g10}_{f00}\GG_{1,0}\!\otimes\!\FF_{0,0}\!\!+\!\!f^{h10}_{f00}\HH_{1,0}\!\otimes\!\FF_{0,0}\!\!+\!\!f^d_{f10}\DD\!\otimes\!\FF_{1,0}\!\!+\!\!f^{d_1}_{f10}\DD_{1}\!\otimes\!\FF_{1,0}\!\!+\!\!f^{d_2}_{f10}\DD_{2}\!\otimes\!\FF_{1,0}.
\end{array}\end{eqnarray}
Applying $\dir_{\bz}$ to $[\DD,\GG_{1,0}]=0,$ using Subclaim \ref{subclaim1} and expression \eqref{E10}, we can simplify $\dir_{\bz}(\GG_{1,0})$ as
\begin{eqnarray}\begin{array}{cccc}\label{G102}
&&g^{e00}_{f10}\EE_{0,0}\!\otimes\!\FF_{1,0}\!+\!g^{e10}_{f00}\EE_{1,0}\!\otimes\!\FF_{0,0}\!+\!g^{f00}_{e10}\FF_{0,0}\!\otimes\!\EE_{1,0}\!+\!g^{f10}_{e00}\FF_{1,0}\!\otimes\!\EE_{0,0} \\[6pt]
&&+g^{g10}_{d}\GG_{1,0}\!\otimes\!\DD+g^{g10}_{d_1}\GG_{1,0}\otimes\DD_{1}+g^{g10}_{d_2}\GG_{1,0}\otimes\DD_{2}+g^{h10}_{d}\HH_{1,0}\!\otimes\!\DD \\[6pt]
&&+g^{h10}_{d_1}\HH_{1,0}\!\otimes\!\DD_{1}+g^{h10}_{d_2}\HH_{1,0}\otimes\!\DD_{2}+g^d_{g10}\DD\!\otimes\GG_{1,0}+g^{d}_{h10}\DD\!\otimes\!\HH_{1,0}\\[6pt]
&&+g^{d_1}_{g10}\DD_{1}\!\otimes\!\GG_{1,0}+g^{d_1}_{h10}\DD_{1}\!\otimes\!\HH_{1,0}+g^{d_2}_{g10}\DD_{2}\!\otimes\!\GG_{1,0}+g^{d_2}_{h10}\DD_{2}\!\otimes\!\HH_{1,0}.
\end{array}\end{eqnarray}
Applying $\dir_{\bz}$ to $[\DD,\HH_{1,0}]=0,$ using Subclaim \ref{subclaim1} and expression \eqref{E10}, we can simplify $\dir_{\bz}(\HH_{1,0})$ as
\begin{eqnarray}\begin{array}{cccc}\label{H102}
&&h^{e00}_{f10}\EE_{0,0}\!\otimes\!\FF_{1,0}\!+\!h^{e10}_{f00}\EE_{1,0}\otimes\FF_{0,0}\!+\!h^{f00}_{e10}\FF_{0,0}\otimes\EE_{1,0}\!+\!h^{f10}_{e00}\FF_{1,0}\otimes\EE_{0,0} \\[6pt]
&&+h^{g10}_{d}\GG_{1,0}\otimes\DD+h^{g10}_{d_1}\GG_{1,0}\otimes\DD_{1}+h^{g10}_{d_2}\GG_{1,0}\otimes\DD_{2}+h^{h10}_d\HH_{1,0}\otimes\DD \\[6pt]
&&+h^{h10}_{d_1}\HH_{1,0}\!\otimes\!\DD_{1}+h^{h10}_{d_2}\HH_{1,0}\otimes\DD_{2}+h^d_{g10}\DD\otimes\GG_{1,0}+h^d_{h10}\DD\otimes\HH_{1,0}\\[6pt]
&&+h^{d_1}_{g10}\DD_{1}\otimes\GG_{1,0}\!+\!h^{d_1}_{h10}\DD_{1}\otimes\HH_{1,0}\!+\!h^{d_2}_{g10}\DD_{2}\otimes\GG_{1,0}\!+\!h^{d_2}_{h10}\DD_{2}\otimes\HH_{1,0}.
\end{array}\end{eqnarray}
Applying $\dir_{\bz}$ to $[\GG_{0,1},\HH_{1,0}]=0,$ we get the following identities \begin{eqnarray}\begin{array}{cccc}\label{gh1}
&&g^{h01}_d=g^{h01}_{d_{1}}=g^{h01}_{d_{2}}=g^d_{h01}=g^{d_{1}}_{h01}=g^{d_{2}}_{h01}=0,\\[6pt]
&&h^{g10}_d=h^{g10}_{d_{1}}=h^{g10}_{d_{2}}=h^d_{g10}=h^{d_{1}}_{g10}=h^{d_{2}}_{g10}=0,\\[6pt]
&&h^{h10}_{d_{2}}=g^{d_{1}}_{g01},\ \ g^{f00}_{e01}=-g^{f01}_{e00}=h^{f10}_{e00}=-h^{f00}_{e10},\\[6pt]
&&h^{d_{2}}_{h10}=g^{g01}_{d_{1}},\ \ g^{e00}_{f01}=-g^{e01}_{f00}=h^{e10}_{f00}=-h^{e00}_{f10}.
\end{array}\end{eqnarray}
Applying $\dir_{\bz}$ to $[\GG_{1,0},\HH_{0,1}]=0,$ we get the following identities
\begin{eqnarray}\begin{array}{cccc}\label{gh2}
&&g^{h10}_d=g^{h10}_{d_{1}}=g^{h10}_{d_{2}}=g^d_{h10}=g^{d_{1}}_{h10}=g^{d_{2}}_{h10}=0,\\[6pt]
&&h^{g01}_{d}=h^{g01}_{d_{1}}=h^{g01}_{d_{2}}=h^d_{g01}=h^{d_{1}}_{g01}=h^{d_{2}}_{g01}=0,\\[6pt]
&&h^{h01}_{d_{1}}=g^{d_{2}}_{g10},\ \ g^{f10}_{e00}=-g^{f00}_{e10}=h^{f00}_{e01}=-h^{f01}_{e00},\\[6pt]
&&h^{d_{1}}_{h01}=g^{g10}_{d_{2}},\ \ g^{e00}_{f10}=-g^{e10}_{f00}=h^{e01}_{f00}=-h^{e00}_{f01}.
\end{array}\end{eqnarray}
Applying $\dir_{\bz}$ to $[\GG_{0,1},\HH_{0,1}]=0$ and $[\GG_{1,0},\HH_{1,0}]=0,$ we have
\begin{eqnarray}\begin{array}{cccc}\label{gh3}
&&h^{f01}_{e00}-h^{f00}_{e01}=g^{f00}_{e01}-g^{f01}_{e00},\ \ h^{h01}_{d_{2}}=g^{d_{2}}_{g01},\ \ \ h^{d_{2}}_{h01}=g^{g01}_{d_{2}},\\[6pt]
&&h^{e00}_{f01}-h^{e01}_{f00}=g^{e01}_{f00}-g^{e00}_{f01},\ \ \,h^{h10}_{d_{1}}=g^{d_{1}}_{g10},\ \ \ h^{d_{1}}_{h10}=g^{g10}_{d_{1}},\\[6pt]
&&h^{f10}_{e00}-h^{f00}_{e10}=g^{f00}_{e10}-g^{f10}_{e00},\ \ h^{e00}_{f10}-h^{e10}_{f00}=g^{e10}_{f00}-g^{e00}_{f10}.
\end{array}\end{eqnarray}
Applying $\dir_{\bz}$ to $[\GG_{0,1},\EE_{0,0}]=\EE_{0,1}$ and $[\HH_{0,1},\EE_{0,0}]=-\EE_{0,1},$ we have
\begin{eqnarray}\begin{array}{llllllll}\label{geh}
e^{e00}_{g01}\!\!\!&=&\!\!\!2g^d_{g01}\!-\!g^{e00}_{f01}=h^{e00}_{f10},\ \ e^{e00}_{h01}\!=\!g^{e00}_{f01}\!+\!2g^d_{h01},\ \ e^{d_{2}}_{e01}\!=\!g^{d_{2}}_{g01}\!-\!g^{d_{2}}_{h01}\!=\!h^{d_{2}}_{h01}\!-\!h^{d_{2}}_{g01},\\[6pt]
e^{h01}_{e00}\!\!\!&=&\!\!\!g^{f01}_{e00}\!+\!2g^{h01}_d\!=\!-\!h^{f01}_{e00}\!-\!2h^{h01}_{d},\ \ e^{e01}_d\!=\!g^{g01}_d\!\!-\!g^{e01}_{f00}\!\!-\!g^{h01}_d\!\!=\!h^{e01}_{f00}\!-\!h^{g01}_d\!+\!h^{h01}_d,\\[6pt]
e^{g01}_{e00}\!\!\!&=&\!\!\!2g^{g01}_d-g^{f01}_{e00}=h^{f01}_{e00}-2h^d_{g01},\ \ \ \ e^{e00}_{h01}=g^{e00}_{f01}+2g^d_{h01}=-h^{e00}_{f01}-2h^d_{h01},\\[6pt]
e^{e01}_{d_{1}}\!\!\!&=&\!\!\!g^{g01}_{d_{1}}-g^{h01}_{d_{1}}\!=\!h^{h01}_{d_{1}}\!-\!h^{g01}_{d_{1}},\ \ \ \ e^d_{e01}\!=\!g^d_{g01}\!-\!g^{f00}_{e01}\!-\!g^d_{h01}\!=\!h^{f00}_{e01}\!-\!h^d_{g01}\!+\!h^d_{h01},\\[6pt]
e^{d_{1}}_{e01}\!\!\!&=&\!\!\!g^{d_{1}}_{g01}-g^{d_{1}}_{h01}=h^{d_{1}}_{h01}-h^{d_{1}}_{g01},\ e^{e01}_{d_{2}}=g^{g01}_{d_{2}}-g^{h01}_{d_{2}}=h^{h01}_{d_{2}}-h^{g01}_{d_{2}}.
\end{array}\end{eqnarray}
Applying $\dir_{\bz}$ to $[\EE_{0,1},\EE_{0,0}]=0,$ we have
\begin{eqnarray}\label{ee}
e^{e00}_{g01}-e^{e00}_{h01}+2e^d_{e01}=0,\ \ \ 2e^{e01}_d+e^{g01}_{e00}-e^{h01}_{e00}=0.
\end{eqnarray}
Combined with equations from \eqref{gh1} to \eqref{ee}, we obtain that
\begin{eqnarray}\begin{array}{llllllll}\label{01}
g^{g01}_{d_{2}}\!\!\!&=&\!\!\!h^{h01}_{d_{2}}=g^{d_{2}}_{g01}=h^{d_{2}}_{h01}=g^{d_{1}}_{g10}=g^{g10}_{d_{1}},\\[6pt]
g^d_{g01}\!\!\!&=&\!\!\!g^{g01}_d\!=\!h^{h01}_d\!=\!h^{h01}_d\!=\!0,\ \ g^{g01}_{d_{1}}\!=\!h^{h01}_{d_{1}}\!=\!g^{d_{1}}_{g01}\!=\!h^{d_{1}}_{h01}\!=\!h^{h10}_{d_{2}}\!=\!h^{d_{2}}_{h10},\\[6pt]
g^{e00}_{f01}\!\!\!&=&\!\!\!-g^{e01}_{f00}=-g^{f00}_{e01}=g^{f01}_{e00}=h^{e01}_{f00}=-h^{e00}_{f01}=h^{f00}_{e01}=-h^{f01}_{e00}\\[6pt]
\!\!\!&=&\!\!\!g^{e00}_{f10}=-g^{e10}_{f00}=-g^{f00}_{e10}=g^{f10}_{e00}=h^{e10}_{f00}=-h^{e00}_{f10}\!=\!h^{f00}_{e10}\!=\!-h^{f10}_{e00}.
\end{array}\end{eqnarray}
By the similar method, applying $\dir_{\bz}$ to $[\GG_{1,0},\EE_{0,0}]=\EE_{1,0}$ and $[\HH_{1,0},\EE_{0,0}]=-\EE_{1,0},$ one has
\begin{eqnarray}\label{10}
g^d_{g10}=g^{g10}_d=h^d_{h10}=h^{h10}_d=0.
\end{eqnarray}
Applying $\dir_{\bz}$ to $[\GG_{1,0},\FF_{0,0}]=-\FF_{1,0}$ and $[\GG_{0,1},\FF_{0,0}]=-\FF_{0,1},$ owing to expression \eqref{F012} and \eqref{F102}, there is
\begin{eqnarray}\begin{array}{llllllll}\label{ff}
&&f^{f00}_{g10}=-g^{f00}_{e10}+2g^d_{g10},\ \ f^{f00}_{h10}=g^{f00}_{e10}+2g^d_{h10},\ \ f^{f10}_{d_1}=-g^{h10}_{d_1}+g^{g10}_{d_1},\\[6pt]
&&f^{h10}_{f00}=g^{e10}_{f00}\,+\,2g^{h10}_d,\ \ f^{g10}_{f00}\!=\!-g^{e10}_{f00}\!+2g^{g10}_{d},\ \ f^{f10}_{d_2}=-g^{h10}_{d_2}+g^{g10}_{d_2},\\[6pt]
&&f^{f00}_{g01}=-g^{f01}_{e00}\!+\!2g^d_{g01},\ \ \,f^{f00}_{h01}\,=g^{f00}_{e01}+2g^d_{h01},\ \ f^{f01}_{d_2}=-g^{h01}_{d_2}+g^{g01}_{d_2},\\[6pt]
&&f^{d_2}_{f01}=g^{d_2}_{g01}\,-g^{d_2}_{h01},\ \ \, f^{g01}_{f00}\!=\!-g^{e01}_{f00}\!+\!2g^{g01}_{d^{0,1}},\ \ \, f^{f01}_{d_1}=-g^{h01}_{d_1}+g^{g01}_{d_1},\\[6pt]
&&f^{h01}_{f00}=g^{e01}_{f00}\,+\,2g^{h01}_d,\ \ \, f^{d_1}_{f01}\,=\, g^{d_1}_{g01}\,-\,g^{d_1}_{h01},\ \ f^{f01}_d\!=\!-g^{f01}_{e00}\!-\!g^{h01}_d\!+\!g^{g01}_d,\\[6pt]
&&f^{d_1}_{f10}=g^{d_1}_{g10}\,-\,g^{d_1}_{h10},\ \ f^{d_2}_{f10}\ =\,g^{d_2}_{g10}\,-\,g^{d_2}_{h10},\ \ \, f^d_{f01}\!=\!-g^{e00}_{f01}\!+\!g^d_{g01}\!-\!g^d_{h01},\\[6pt]
&&f^d_{f10}=-g^{e00}_{f10}+g^d_{g10}-g^d_{h10},\ \ \ \ \ \ f^{f10}_d=-g^{f10}_{e00}-g^{h10}_d+g^{g10}_d.
\end{array}\end{eqnarray}

Redenoting
$\dir_{\bz}+u_{inn}$ by $\dir_{\bz}$, where $u=g^{g01}_{d_{2}}\DD_{2}\otimes\DD_{2}+g_{d_{1}}^{g01}\DD_{2}\otimes\DD_{1}+g^{d_{1}}_{g01}\DD_{1}\otimes\DD_{2}+g^{d_{1}}_{g10}\DD_{1}\otimes\DD_{1}$, then using equation \eqref{gh1}--\eqref{ff}, one has
\begin{eqnarray*}
\dir_{\bz}(\GG_{0,1})\!\!\!&=&\!\!\! g^{f00}_{e01}\big(\FF_{0,0}\otimes\EE_{0,1}\!-\!\FF_{0,1}\otimes\EE_{0,0}\!-\!\EE_{0,0}\otimes\FF_{0,1}+\EE_{0,1}\otimes\FF_{0,0}\big),\\
\dir_{\bz}(\HH_{0,1})\!\!\!&=&\!\!\!-g^{f00}_{e01}\big(\FF_{0,0}\otimes\EE_{0,1}\!-\!\FF_{0,1}\otimes\EE_{0,0}-\EE_{0,0}\otimes\FF_{0,1}+\EE_{0,1}\otimes\FF_{0,0}\big),\\
\dir_{\bz}(\EE_{0,1})\!\!\!&=&\!\!\! g^{f00}_{e01}\big(\EE_{0,0}\otimes\GG_{0,1}\!-\!\EE_{0,0}\otimes\HH_{0,1}\!-\!\EE_{0,1}\otimes\DD+\GG_{0,1}\otimes\EE_{0,0}\!-\!\HH_{0,1}\otimes\EE_{0,0}\!-\!\DD\otimes\EE_{0,1}\big),\\
\dir_{\bz}(\FF_{0,1})\!\!\!\!&=&\!\!-g^{f00}_{e01}\big(\FF_{0,0}\otimes\GG_{0,1}\!-\!\FF_{0,0}\otimes\HH_{0,1}\!-\!\FF_{0,1}\otimes\DD+\GG_{0,1}\otimes\FF_{0,0}\!-\!\HH_{0,1}\otimes\FF_{0,0}-\DD\otimes\FF_{0,1}\big),\\
\dir_{\bz}(\GG_{1,0})\!\!\!&=&\!\!\! g^{f00}_{e01}\big(\FF_{0,0}\otimes\EE_{1,0}\!-\!\FF_{1,0}\otimes\EE_{0,0}\!-\!\EE_{0,0}\otimes\FF_{1,0}+\EE_{1,0}\otimes\FF_{0,0}\big),\\
\dir_{\bz}(\HH_{1,0})\!\!\!&=&\!\!\!-g^{f00}_{e01}\big(\FF_{0,0}\otimes\EE_{1,0}\!-\!\FF_{1,0}\otimes\EE_{0,0}\!-\!\EE_{0,0}\otimes\FF_{1,0}+\EE_{1,0}\otimes\FF_{0,0}\big),\\
\dir_{\bz}(\EE_{1,0})\!\!\!&=&\!\!\! g^{f00}_{e01}\big(\EE_{0,0}\otimes\GG_{1,0}\!-\!\EE_{0,0}\otimes\HH_{1,0}\!-\!\EE_{1,0}\otimes\DD+\GG_{1,0}\otimes\EE_{0,0}\!-\!\HH_{1,0}\otimes\EE_{0,0}\!-\!\DD\otimes\EE_{1,0}\big),\\
\dir_{\bz}(\FF_{1,0})\!\!\!&=&\!\!\!-g^{f00}_{e01}\big(\FF_{0,0}\otimes\GG_{1,0}\!-\!\FF_{0,0}\otimes\HH_{1,0}\!-\!\FF_{1,0}\otimes\DD+\GG_{1,0}\otimes\FF_{0,0}\!-\!\HH_{1,0}\otimes\FF_{0,0}-\DD\otimes\FF_{1,0}\big).\\
\end{eqnarray*}
By careful observations and patient calculations, taking $\tilde{\dir}=g^{f00}_{e01}(\EE_{\bz}\otimes\FF_{\bz}+\FF_{\bz}\otimes\EE_{\bz}+1/2\DD\otimes\DD)$, one can check that
\begin{eqnarray*}
\tilde{\dir}(\EE_{0,0})\!\!\!&=&\!\!\! g^{f00}_{e01}(\EE_{\bz}\otimes\DD+\DD\otimes\EE_{\bz}-1/2(2\EE_{\bz}\otimes\DD+2\DD\otimes\EE_{\bz}))=0,\\
\tilde{\dir}(\FF_{0,0})\!\!\!&=&\!\!\!g^{f00}_{e01}(-\DD\otimes\EE_{\bz}\!-\!\EE_{\bz}\otimes\DD\!+1/2(2\EE_{\bz}\otimes\DD\!+2\DD\otimes\EE_{\bz}))=0,\\
\tilde{\dir}(\DD)\!\!\!&=&\!\!\!0,\\
\tilde{\dir}(\GG_{0,1})\!\!\!&=&\!\!\!\dir_{\bz}(\GG_{0,1}),\ \ \tilde{\dir}(\HH_{0,1})=\dir_{\bz}(\HH_{0,1}),\ \ \tilde{\dir}(\EE_{0,1})=\dir_{\bz}(\EE_{0,1}),\ \ \tilde{\dir}(\FF_{0,1})=\dir_{\bz}(\FF_{0,1}),\\
\tilde{\dir}(\GG_{1,0})\!\!\!&=&\!\!\!\dir_{\bz}(\GG_{1,0}),\ \ \tilde{\dir}(\HH_{1,0})=\dir_{\bz}(\HH_{1,0}),\ \ \tilde{\dir}(\EE_{1,0})=\dir_{\bz}(\EE_{1,0}),\ \ \tilde{\dir}(\FF_{1,0})=\dir_{\bz}(\FF_{1,0}).
\end{eqnarray*}
Since $\WL_{0,1}={\rm Span}_\C\{\GG_{0,1},\HH_{0,1},\EE_{0,1},\FF_{0,1}\}$ and $\WL_{1,0}={\rm Span}_\C\{\GG_{1,0},\HH_{1,0},\EE_{1,0},\FF_{1,0}\}$, we complete the proof of Subclaim \ref{subclaim2}.\QED\vskip5pt
\begin{subclai}
\label{subclaim3} $\dir_{\bz}(\WL)=0.$
\end{subclai}
\noindent{\it Proof}. According to the fact that the algebra $\WL$ is generated by the set $$\{\EE_{0,0},\,\FF_{0,0},\,\DD_{1},\,\DD_{2},\, \EE_{1,0},\,\FF_{1,0},\,\EE_{0,1},\,\FF_{0,1},\,\GG_{0,n},\,\GG_{n,0}\,|\,n\in\Z_{>0}\},$$ and using all of the above Subclaims, we only need to prove $$\dir_{\bz}(\GG_{0,n})=\dir_{\bz}(\GG_{n,0})=0\ \ {\rm for\ all}\ n\in\Z_{>1}.$$

For some $n\geq2$. Applying $\dir_{\bz}$ to $[\DD,\GG_{0,n}]=0,$ using Subclaim \ref{subclaim1}, we can simplify $\dir_{\bz}(\GG_{0,n})$ as
\begin{eqnarray}\label{G012}
&&\ \ \alpha_m\EE_{0,m}\otimes\FF_{0,n-m}+\alpha_m^{+}\FF_{0,m}\otimes\EE_{0,n-m}+\beta_m\GG_{0,m}\otimes\GG_{0,n-m}\nonumber\\
&&+\beta_m^{+}\HH_{0,m}\otimes\HH_{0,n-m}+\gamma_m\GG_{0,m}\otimes\HH_{0,n-m}+\gamma_m^+\HH_{0,m}\otimes\GG_{0,n-m}\nonumber\\
&&+\xi\GG_{0,n}\otimes\DD+\xi_1\GG_{0,n}\otimes\DD_{1}+\xi_2\GG_{0,n}\otimes\DD_{2}+\zeta\HH_{0,n}\otimes\DD\\
&&+\zeta_1\HH_{0,n}\otimes\DD_{1}+\zeta_2\HH_{0,n}\otimes\DD_{2}+\xi^+\DD\otimes\GG_{0,n}+\xi_1^+\DD_{1}\otimes\GG_{0,n}\nonumber\\
&&+\xi_2^+\DD_{2}\otimes\GG_{0,n}+\zeta^+\DD\otimes\HH_{0,n}+\zeta_1^+\DD_{1}\otimes\HH_{0,n}+\zeta_2^+\DD_{2}\otimes\HH_{0,n},\nonumber \end{eqnarray}
where, $m\in\Z,0\leq m\leq n$ and $\beta_0=\beta_n=\beta_0^+=\beta_n^+=\gamma_0=\gamma_n=\gamma_0^+=\gamma_n^+=0$.  Without confusion, we would assume $\alpha_{n+1}=\alpha_{n+1}^+=0$ for convenience.

Applying $\dir_{\bz}$ to $[\GG_{0,1},\GG_{0,n}]=0$, there is $\xi_2=\zeta_2=\xi_2^+=\zeta_2^+$ and
\begin{eqnarray*}
\alpha_{m+1}=\alpha_m,\ \ \alpha_{m+1}^+=\alpha_m^+,
\end{eqnarray*}
where $0\leq m\leq n$. Thus, $\alpha_{m}=\alpha_{m}^+=0$ for $0\leq m\leq n$.

Applying $\dir_{\bz}$ to $[\HH_{1,0},\GG_{0,n}]=0$, we have
\begin{eqnarray*}
\beta_{m}^+=\gamma_m=\gamma_m^+=\zeta=\zeta_1=\zeta^+=\zeta_1^+=\xi_1=\xi_1^+=0,
\end{eqnarray*}
where $0< m< n$.

Now, we can simplify $\dir_{\bz}(\GG_{0,n})$ as
\begin{eqnarray}\label{G0121} \beta_m\GG_{0,m}\otimes\GG_{0,n-m}+\xi\GG_{0,n}\otimes\DD+\xi^+\DD\otimes\GG_{0,n},
\end{eqnarray}
where $m\in\Z,0\leq m\leq n$ and $\beta_0=\beta_n=0$.

Since $\dir_{\bz}(\EE_{1,0})=\dir_{\bz}(\GG_{1,0})=\dir_{\bz}(\GG_{0,1})=0$, then $\dir_{\bz}(\EE_{m,n})=0$. Applying $\dir_{\bz}$ to $[\EE_{0,0},\GG_{0,n}]=\EE_{0,n}$, we have $\EE_{0,0}\cdot\dir_{\bz}(\GG_{0,n})=0$, which equal to
\begin{eqnarray}
\beta_m=\xi=\xi^+=0,
\end{eqnarray}
for $0<m<n$. Thus, $\dir_{\bz}(\GG_{0,n})=0$.

Similarly, by the same method, we can obtain $\dir_{\bz}(\GG_{n,0})=0$.

 Then this subclaim follows.\QED\vskip5pt

By now, we have completed the proof of Claim \ref{clai3}.\QED\vskip5pt
\begin{clai}\label{clai4}
For any $\dir\in{\rm Der}(\WL,\VV)$, \eqref{summable} is a finite sum.
\end{clai}
\noindent{\it Proof}.\ \ Since $\dir=\mbox{$\sum\limits_{\bk\in\Z^2}\dir_\bk$},$ by the above claims, one can suppose $\dir_{m_{1},m_{2}}=(v_{m_{1},m_{2}})_{\rm
inn}$ for some $v_{m_{1},m_{2}}\in\VV_{m_{1},m_{2}}$ and $(m_{1},m_{2})\in\Z^{2}_{\geq 0}.$ If $\Gamma'=\{(m_{1},m_{2})\in\Z^{2}_{> 0}|\,v_{m_{1},m_{2}}\neq0\}$ is an infinite set, by linear algebra, there exists $\rho\in\texttt{T}$ such that $\rho(m_{1},m_{2})\neq0$ for $(m_{1},m_{2})\in \Gamma',$ Then $\dir(\rho)=\sum_{(m_{1},m_{2})\in\Gamma'}\rho(m_{1},m_{2})v_{m_{1},m_{2}}$ is an infinite sum, which is not an element in $\VV.$ This is a contraction with the fact that $\dir\in{\rm Der}(\WL,\VV).$

This proves Claim \ref{clai4} and Proposition \ref{proposition}.\QED\vskip5pt

\begin{lemm}\label{lemma5v} Suppose $v\in\VV$ such that
$x\cdot v\in {\rm Im}(1-\tau)$ for all $x\in\WL.$ Then $v\in {\rm
Im}(1-\tau)$.\end{lemm}
\noindent{\it Proof}.\ \ It is easy to see that ${\rm Im}(1-\tau)={\rm ker}(1+\tau)$. Then for any $v\in\VV$ such that $\WL\cdot v\in {\rm Im}(1-\tau)$, one has $(1+\tau)(\WL\cdot v)=0$. Noting that $\tau$ commutes with the action of $\WL$ on $\VV$, we obtain $\WL\cdot (1+\tau)v=(1+\tau)(\WL\cdot v)=0$, which together with Lemma \ref{Legr}, forces $(1+\tau)v=0$. Then this lemma follows from the fact that ${\rm Ker}(1+\tau)={\rm Im}(1-\tau)$.\QED\vskip5pt

\vskip10pt \ni{\it Proof of Theorem \ref{main}}\ \ Let $(\WL
,[\cdot,\cdot],\dir)$ be a Lie bialgebra structure on $\WL$. By
(\ref{bLie-d}), (\ref{deriv}) and Proposition \ref{proposition},
$\D=\D_r$ is defined by (\ref{D-r}) for some $r\in\WL\otimes\WL.$ By
(\ref{cLie-s-s}), ${\rm Im}\,\D\subset{\rm Im}(1-\tau).$ Thus by
Lemma \ref{lemma5v}, $r\in{\rm Im}(1-\tau).$ Then (\ref{cLie-s-s}),
(\ref{add-c}) and Corollary \ref{coro1} show that $c(r)=0.$ Then
Definition \ref{def2} says that $(\WL ,[\cdot,\cdot],\D)$ is a
triangular coboundary Lie bialgebra.\QED\vskip7pt

\noindent{\bf Acknowledgements}\ \ The authors would sincerely like to thank the referee for the invaluable comments, in particular providing the much simpler proofs of Subclaim 1 in Claim 3 of Proposition 2.4 and Lemma 2.5, which help us avoid the heavy computations.

\vskip8pt

\begin{thebibliography}{9999}\vskip0pt\small
\parindent=2ex\parskip=-1pt\baselineskip=-1pt

\bibitem{B} S. Berman, Y. Gao and Y.S. Krylyuk, Quantum tori and the strcutre of elliptic quasi-simple Lie algebras, {\it J. Funct. Anal.} {\bf135} (1996), 339--389.

\bibitem{BJ} S. Berman and J. Szmigielski, Principal realization for extended affine Lie algebra of type $\frak{sl_2}$ with coordinates in a simple quantum torus with two generators, arXiv:hep-th/9805016v1.

\bibitem{CS} Y.Cheng,Y.Shi, Lie Bialgebra Structures on the q-Analog Virasoro-Like Algebras, {\it Communications in Algebra} {\bf 37}(4) (2009), 1264--1274.

\bibitem{D1} V.G. Drinfeld, Constant quasiclassical solutions of the
Yang-Baxter quantum equation, {\it Soviet Math. Dokl.} {\bf28}(3) (1983).
667--671.

\bibitem{D2} V.G. Drinfeld, Quantum groups, in: {\it Proceeding of the
International Congress of Mathematicians}, Vol.~1, 2, Berkeley,
Calif.~1986, Amer.~Math.~Soc., Providence, RI, 1987, pp.~798--820.

\bibitem{G} C. Grunspan, Quantizations of the Witt algebra and of simple Lie
algebras in characteristic $p$, {\it J. Alg.} {\bf 280} (2004).
145--161.

\bibitem{G1} Y. Gao, Vertex operators arising from the homogeneous realization for d, {\it Comm. Math.
Phys.} {\bf211} (2000), 745--777.

\bibitem{G2} Y. Gao, Representation of extended affine Lie algebras coordinatized by certain quantum tori,
{\it Composito Mathematica} {\bf123} (2000), 1--25.

\bibitem{J} H.P. Jakobson and V.G. Kac, A new class of unitarizable highest weight representations of
infinite-dimensional Lie algebras II, {\it J. Funct. Anal.} {\bf82} (1989), 69--90.

\bibitem{LSX} J. Li, Y. Su, B. Xin, Lie bialgebras of a family of Block
type, {\it Chinese Annals of Math. (Series.B)} {\bf 29} (2008),
487--500.

\bibitem{M} Kei Miki, Integrable irreducible higeest weight modules for $\frak{sl}_2(C_p[x^{\pm1}, y^{\pm1}])$, {\it Osaka J. Math.} {\bf41} (2004), 295--326

\bibitem{MY} Y.I. Manin, Topics in noncommutative geometry, {\it Princeton Univ.} Press, 1991.

\bibitem{Mi} W. Michaelis, A class of infinite-dimensional Lie
bialgebras containing the Virasoro algebras, {\it Adv. Math.}
{\bf107} (1994), 365--392.

\bibitem{NT} S.H. Ng, E.J.~Taft, Classification of the Lie bialgebra
structures on the Witt and Virasoro algebras, {\it J. Pure Appl.
Alg.} {\bf151} (2000), 67--88.

\bibitem{SS} G. Song, Y. Su, Lie bialgebras of
generalized Witt type, {\it Science in China: Series A Mathematics.}
{\bf49}(4) (2006), 533--544.

\bibitem{T} E.J. Taft, Witt and Virasoro algebras as Lie bialgebras, {\it
J. Pure Appl.~Algebra} {\bf87} (1993), 301--312.

\bibitem{V} M. Varagnolo and E. Vasserot, Double loop algebras and Fock space, {\it Invent. Math.} {\bf133} (1998), 133--159.

\bibitem{WSS} Y. Wu, G. Song, Y. Su, Lie bialgebras of generalized
Virasoro-like type, {\it Acta Mathematica Sinica, English Series}
{\bf22} (2006), 1915--1922.



\end{thebibliography}
\end{document}